%
%


\input amstex

\magnification 1200
\loadmsbm
\parindent 0 cm

\define\nl{\bigskip\item{}}
\define\snl{\smallskip\item{}}
\define\inspr #1{\parindent=20pt\bigskip\bf\item{#1}}
\define\iinspr #1{\parindent=27pt\bigskip\bf\item{#1}}
\define\einspr{\parindent=0cm\bigskip}

\define\ot{\otimes}

\define\tr{\triangleright}
\define\tl{\triangleleft}
\define\btr{\blacktriangleright}
\define\btl{\blacktriangleleft}

\input amssym
\input amssym.def


\centerline{\bf Algebraic quantum hypergroups II.}
\centerline{\bf Constructions and examples.}
\bigskip\bigskip
\centerline{\it L.\ Delvaux\,$^{\text{1}}$ and A.\ Van Daele\,$^{\text{2}}$} 
\bigskip\bigskip\bigskip
{\bf Abstract} 
\bigskip
Let $G$ be a group and let $A$ be the algebra of complex functions on $G$ with finite support. The product in $G$ gives rise to a coproduct $\Delta$ on $A$ making the pair $(A,\Delta)$ a multiplier Hopf algebra. In fact, because there exist integrals, we get an algebraic quantum group as introduced and studied in [VD2].
\snl
Now let $H$ be a finite subgroup of $G$ and consider the subalgebra $A_1$ of functions in $A$ that are constant on double cosets of $H$. The coproduct in general will not leave this algebra invariant but we can modify $\Delta$ and define $\Delta_1$ as
$$\Delta_1(f)(p,q)=\frac{1}{n}\sum_{r\in H} f(prq)$$
where $f\in A_1$, $p,q\in G$ and where $n$ is the number of elements in the subgroup $H$. Then $\Delta_1$ will leave the subalgebra invariant (in the sense that the image is in the multiplier algebra $M(A_1\ot A_1)$ of the tensor product). However,  it will no longer be an algebra map. So, in general we do not have an algebraic quantum group but a so-called {\it algebraic quantum hypergroup} as introduced and studied in [De-VD].
\snl
Group-like projections in a $^*$-algebraic quantum group $A$ (as defined and studied in [L-VD]) give rise, in a natural way, to $^*$-algebraic quantum hypergroups, very much like subgroups do as above for a $^*$-algebraic quantum group associated to a group (again see [L-VD]). In this paper we push these results further. On the one hand, we no longer assume the $^*$-structure as in [L-VD] while on the other hand, we allow the group-like projection to belong to the multiplier algebra $M(A)$ of $A$ and not only to $A$ itself.  Doing so, we not only get some well-known earlier examples of algebraic quantum hypergroups but also some interesting new ones.
\nl
\nl
January 2010 ({\it Version} 1.0)
\vskip 1 cm
\hrule
\bigskip\parindent 0.3 cm
\item{$^{\text{1}}$} Department of Mathematics, Hasselt University, Agoralaan, B-3590 Diepenbeek (Belgium). {\it E-mail}: Lydia.Delvaux\@uhasselt.be
\item{$^{\text{2}}$} Department of Mathematics, K.U.\ Leuven, Celestijnenlaan 200B (bus 2400),
B-3001 Heverlee (Belgium). {\it E-mail}: Alfons.VanDaele\@wis.kuleuven.be
\parindent 0 cm
 
\newpage


\bf 0. Introduction \rm
\nl
Let us start with recalling the basic notions used in this paper.
\snl
A {\it multiplier Hopf algebra} is a pair $(A,\Delta)$ of an algebra $A$ and a coproduct $\Delta$ on $A$ satisfying certain properties. We assume that the algebra is an algebra over the field $\Bbb C$ of complex numbers. It need not be unital, but the product is always assumed to be non-degenerate (as a bilinear map). Then we can consider the multiplier algebra $M(A)$. It contains $A$ as an essential two-sided ideal. The coproduct is a homomorphism $\Delta$ from $A$ to $M(A\ot A)$, the multiplier algebra of the tensor product. If $A$ and $B$ are two algebras with a non-degenerate product and $\gamma:A\to M(B)$ is a homomorphism, it is called non-degenerate if $\gamma(A)B=B\gamma(A)=B$. Then $\gamma$ has a unique extension to a unital homomorphism, still denoted by $\gamma$, from $M(A)$ to $M(B)$. This result can be applied in several situations. One assumes e.g.\ that the coproduct $\Delta$ is non-degenerate. Then one also has that the maps $\Delta\ot \iota$ and $\iota\ot \Delta$ from $A\ot A$ to $M(A\ot A\ot A)$ (where $\iota$ is used for the identity map) are non-degenerate. Using the extensions, one can formulate coassociativity as $(\Delta\ot \iota)\Delta=(\iota\ot\Delta)\Delta$. It is also assumed that the linear maps $a\ot b\mapsto \Delta(a)(1\ot b)$ and $a\ot b\mapsto (a\ot 1)\Delta(b)$ are injective, map into $A\ot A$ and have all of $A\ot A$ as range. Here, $1$ is used to denote the identity in the multiplier algebra $M(A)$ of $A$. One can show the existence of a counit and an antipode. We use $S$ to denote the antipode and $\varepsilon$ for the counit. A multiplier Hopf algebra is regular if the antipode $S$ is a bijective map from $A$ to itself. If $A$ is a $^*$-algebra, we assume that $\Delta$ is a $^*$-homomorphism and then the pair $(A,\Delta)$ is called a multiplier Hopf $^*$-algebra. In this case regularity is automatic.
\snl
For the definition and the basic theory on multiplier Hopf algebras, we refer to [VD1]. 
\snl
An {\it algebraic quantum group} is a regular multiplier Hopf algebra with integrals. We use $\varphi$ for a left integral and $\psi$ for a right integral. They are faithful linear functionals on $A$ satisfying left, respectively right invariance. Left invariance of $\varphi$ is expressed as $(\iota\ot\varphi)\Delta(a)=\varphi(a)1$ in the multiplier algebra $M(A)$ for all $a\in A$. Similarly for the right invariance of $\psi$. Integrals are unique (up to a scalar). There exists a multiplier $\delta$, called the modular element, defined and characterized by $(\varphi\ot\iota)\Delta(a)=\varphi(a)\delta$ for all $a\in A$. It also satisfies $(\iota\ot\psi)\Delta(a)=\psi(a)\delta^{-1}$ as well as $\varphi(S(a))=\varphi(a\delta)$ for all $a$. Finally, there exist the so-called modular automorphisms $\sigma$ and $\sigma'$ of $A$, defined by $\varphi(ab)=\varphi(b\sigma(a))$ and $\psi(ab)=\psi(b\sigma'(a))$ for all $a,b\in A$. 
\snl
The {\it dual} $\widehat A$ of an algebraic quantum group is defined as the  subspace of linear functionals on $A$ of the form $\varphi(\,\cdot\,a)$ where $a$ runs through $A$. It can be given the structure of an algebraic quantum group. The product in $\widehat A$ is defined dual to the coproduct of $A$ and the coproduct $\widehat\Delta$ on $\widehat A$ is defined dual to the product in $A$. The pair $(\widehat A,\widehat\Delta)$  is called the dual of $(A,\Delta)$ (in the sense of algebraic quantum groups). We will occasionally use $\langle a,\omega\rangle$ for the value $\omega(a)$ when $a\in A$ and $\omega\in \widehat A$. The pair $(A,\widehat A)$ is a dual pair of multiplier Hopf algebras (as defined and studied in [Dr-VD]).
\snl
There are natural left and right actions of $\widehat A$ on $A$ defined by
$$\omega\btr a=(\iota\ot\omega)\Delta(a)\qquad\qquad\text{and}\qquad\qquad a\btl\omega=(\omega\ot\iota)\Delta(a)$$
where $a\in A$ and $\omega \in \widehat A$. These actions make $A$ into a unital $\widehat A$-bimodule. The actions extend in a natural way to the multiplier algebra $M(\widehat A)$ because the actions are unital (see e.g.\ Proposition 3.3 in [Dr-VD-Z]). Similarly, $\widehat A$ is a unital $A$-bimodule with these properties.  
\snl
One can extend the pairing to $A\times M(\widehat A)$ and to $M(A)\times \widehat A$. This can be done in a natural way using the fact that the four modules as considered are unital. So, one can pair elements of $A$ with elements of $M(\widehat A)$ and elements of $M(A)$ with elements of $\widehat A$. However, in general, it is not possible to pair elements of $M(A)$ with elements of $M(\widehat A)$. There are a few exceptions though. It makes sense to consider e.g.\ $\langle m,1 \rangle$ and $\langle m,\delta_{\widehat A} \rangle$ for $m\in M(A)$. The reason is that $1$ and $\delta_{\widehat A}$ are group-like elements in the sense that they give a homomorphism on $A$, namely $\varepsilon$ in the first case and $\varepsilon\circ\sigma$ in the second case. 
\snl
If $(A,\Delta)$ is an algebraic quantum group, $A$ a $^*$-algebra and $\Delta$ a $^*$-homomorphism, then it is called a $^*$-algebraic quantum group. The dual $(\widehat A,\widehat \Delta)$ is again a $^*$-algebraic quantum group. The involution in $\widehat A$ is defined by the formula $\omega^*(a)=\omega(S(a)^*)^-$ where $^-$ denotes complex conjugation in $\Bbb C$. Often in this case, it is assumed that the integrals are positive linear functionals. A linear functional $\omega$ is called positive if $\omega(a^*a)\geq 0$ for all $a$ in the algebra. It can be shown that when a positive left integral exists, also a positive right integral exists and that the integrals on the dual can also be chosen to be positive (see e.g.\ [DC-VD]).
\snl
We refer to [VD2] for the theory of these algebraic quantum groups.
\snl
An {\it algebraic quantum hypergroup} is like an algebraic quantum group, but without the assumption that the coproduct is an algebra map. More precisely, it is a pair $(A,\Delta)$ where $A$ is an algebra as before and where $\Delta$ is a coassociative linear map from $A$ to $M(A\ot A)$ satisfying some other properties as in the case of multiplier Hopf algebras (see Section 1 in [De-VD]). First of all, it is assumed that there is a counit $\varepsilon$, as well as a faithful left integral $\varphi$. These objects are defined just as in the case of multiplier Hopf algebras. An antipode, relative to $\varphi$ is a bijective linear map $S:A\to A$ satisfying 
$$S((\iota\ot\varphi)(\Delta(a)(1\ot b)))=(\iota\ot\varphi)((1\ot a)\Delta(b))$$
for all $a,b\in A$. The map is uniquely defined by this equation. It is assumed to be a anti-homomorphism. Then the pair $(A,\Delta)$ is an algebraic quantum hypergroup. Also the dual $(\widehat A,\widehat\Delta)$ of $(A,\Delta)$ can be constructed as in the case of algebraic quantum groups. It is again an algebraic quantum hypergroup.
\snl
Also here we can consider the involutive case and we speak about $^*$-algebraic quantum hypergroups.
\snl
We refer to [De-VD] for the theory of algebraic quantum hypergroups.
\nl
Now let $(A,\Delta)$ be an algebraic quantum group.
\snl
A {\it group-like projection} in $A$ is a non-zero idempotent element $h\in A$ satisfying $S(h)=h$ (where $S$ is the antipode of $A$) and 
$$\Delta(h)(1\ot h)=h\ot h \qquad\qquad \text{and} \qquad\qquad \Delta(h)(h\ot 1)=h\ot h.$$
The concept has been introduced in [L-VD] in the context of $^*$-algebraic quantum groups. In this paper we do not require the existence of an involution and therefore, the assumptions are a little different. We also consider the more general situation where $h$ is allowed to belong to the multiplier algebra $M(A)$ and not necessarily to $A$ itself (as in [L-VD]).
\snl
A simple example comes from a pair $(G,H)$ of a group $G$ with a  finite subgroup $H$. Here $A$ is the algebra of complex functions with finite support in $G$ and $\Delta$ is defined by $\Delta(f)(p,q)=f(pq)$ where $f\in A$ and $p,q\in G$. If we let $h$ be the function on $G$ that is $1$ on elements of $H$ and $0$ everywhere else, we get a group-like projection. The element $h$ is in $A$ because $H$ is assumed to be a finite and it satisfies the requirements because $H$ is a subgroup. 
\snl
If the subgroup is no longer assumed to be finite, one can still define $h$ as above. Now it will no longer belong to $A$ but rather it will be an element in the multiplier algebra $M(A)$. On the other hand, the main defining conditions will still be satisfied. To give a meaning to them, one first has to extend the coproduct and the antipode to the multiplier algebra, but this can be done in a unique way by a general procedure as we mentioned already.
\snl
Suppose for a moment that we have a group-like projection $h$ in $A$. Then, there are two ways to associate an algebraic quantum hypergroup.
\snl
First, one takes the subalgebra $A_0$ of $A$, defined as $hAh$ . This is an algebra with identity (namely $h$ because we assume here that $h\in A$). The new coproduct $\Delta_0$ is obtained by cutting down the original one as
$$\Delta_0(a)=(h\ot h)\Delta(a)(h\ot h)$$
for $a\in A_0$. Now we have $\Delta_0:A_0\to A_0\ot A_0$, the map is coassociative and linear, but it is no longer an algebra map (except when e.g.\ $h$ is central in $A$). The pair $(A_0,\Delta_0)$ is a an algebraic quantum hypergroup of compact type (see  Definition 1.13 in [De-VD] and Theorem 2.7 in [L-VD]). 
\snl
The second case is dual to the first one. Consider the dual $\widehat A$ of $A$ (in the sense of algebraic quantum groups) and let $k=\varphi(\,\cdot\,h)$ where $\varphi$ is a left integral on $A$, normalized so that $\varphi(h)=1$. Then $k$ is a group-like projection in the dual $\widehat A$. One can take the associated algebraic quantum hypergroup of compact type as above and then the dual (in the sense of algebraic quantum hypergroups). This will give an algebraic quantum hypergroup of discrete type. The underlying algebra is the subalgebra $A_1$ of $A$ defined as the intersection of the left and the right legs of $\Delta(h)$. Also here it is needed  to cut down the coproduct in a certain way. We  refer to Section 3 in [L-VD]; see also Section 3 in [De-VD].
\snl
If however we take a group-like idempotent $h$ in the multiplier algebra $M(A)$ and not in $A$ itself, we can still consider the first construction but not the second one. The Fourier transform $k$ will no longer be defined in $\widehat A$ or even in $M(\widehat A)$. This takes us to the content of the paper.
\nl
\it Content of the paper \rm
\nl
In {\it Section} 1 of this paper, we will see how we can get around the problem mentioned above. We start with a group-like projection $k$ in $M(\widehat A)$. This gives rise to commuting conditional expectations $E$ and $E'$ on $A$. The intersection of the ranges $E(A)$ and $E'(A)$ of these maps turns out to be a subalgebra $A_1$ and by cutting down the coproduct $\Delta$ we get a coproduct $\Delta_1$ on this subalgebra making it into an algebraic quantum hypergroup. We show that this procedure extends the second construction above in the case where the group-like idempotent $k$ is in the algebra $\widehat A$ (and hence is the Fourier transform of a group-like idempotent $h\in A$). Still, also in the more general case, the two constructions yield dual algebraic quantum hypergroups. The main results are found in Theorem 1.15  and Theorem 1.16.
\snl
Let us mention here that in [K] one can also find a construction of quantum hypergroups using conditional expectations on compact quantum groups. In that paper, the author describes non-trivial examples of finite-dimensional quantum hypergroups, coming from Kac algebras.
\snl
In {\it Section} 2, we treat some special cases and we discuss examples to illustrate the results of Section 1. We start with the simplest (motivating) case of a finite subgroup of a group as mentioned already in the abstract and earlier in this introduction. Some other of our examples are still group related. Others work with more general Hopf algebras or multiplier Hopf algebras. In a few examples, the group-like projection $k$ does not belong to the algebra but only to the multiplier algebra. This gives new cases, not treated in [L-VD]. As a special case, we consider the situation with a central idempotent. Then the resulting algebraic quantum hypergroups are actually algebraic quantum groups.

\nl\nl
\bf Acknowledgements   \rm
\nl
We are indebted to Kenny De Commer for letting us use his example to illustrate the results in our paper (cf.\ Example 2.7 in Section 2).
\snl
\nl


\bf 1. Constructions of algebraic quantum hypergroups \rm
\nl
Throughout the paper, $(A,\Delta)$ will be a ($^\ast$-)algebraic quantum group. We will use the notations for the various objects related to this algebraic quantum group as reviewed in the introduction. We refer to the basic references for multiplier Hopf algebras and algebraic quantum groups,  also mentioned in the introduction.
\nl
\it Group-like idempotents in $M(A)$ \rm
\nl
The following definition is basically introduced in [L-VD, Definition 1.1]. Remark that in that paper, one only works with $^*$-algebraic quantum groups (with positive integrals) whereas here, we consider also the non-involutive case. Moreover now we allow elements in the multiplier algebra. That is the reason why the following definition is slightly different from the one in [L-VD].

\inspr{1.1} Definition \rm
A non-zero idempotent $h \in M(A)$ is called  {\it group-like} if $S(h) = h$ and
$$\Delta(h)(1 \otimes h) = h \otimes h \qquad\qquad\text{and}\qquad\qquad \Delta(h) (h \otimes 1) = h\otimes h.$$
In the $^\ast$-case, we furthermore require that $h^\ast = h$.
\einspr

Because we assume that $h$ is invariant under the antipode $S$ and because $S$ is a anti-homomorphisms that flips the coproduct, it is very easy to see that also
$$(1 \otimes h) \Delta(h) = h \otimes h \qquad\qquad\text{and}\qquad\qquad  (h \otimes 1) \Delta(h)=h\ot h.$$
And applying the counit on any of these equations, we get $\varepsilon(h) = 1$ (since $h$ is assumed to be non-zero).
\snl
In [L-VD] it is shown that in the case of a self-adjoint group-like idempotent in a $^*$-algebraic quantum group, it is sufficient to assume only that e.g.\ $\Delta(h)(1 \otimes h) = h \otimes h$. The other three equalities will all follow, as well as the fact that $S(h)=h$. See Proposition 1.6 in [L-VD]. 
\snl
In Proposition 1.7 of [L-VD], it is also shown that $\sigma(h)=h$. For this however, it seems essential that we have the involutive structure and positive integrals. In our situation, we find the following results.

\inspr{1.2} Proposition \rm 
Let $h$ be a group-like idempotent in $M(A)$. If $\varepsilon(\sigma(h))\neq 0$ then $\sigma(h)=h$ (and then $\varepsilon(\sigma(h))=1$). If $\varepsilon(\sigma(h))=0$, then $h\sigma(h)=0$ and then $\varphi(hah)=0$ for all $a\in A$.

\snl\bf Proof\rm:
We know that $\Delta(\sigma(h))=(S^2\ot \sigma)\Delta(h)$. Therefore we get
$$\Delta(\sigma(h))(1\ot \sigma(h))=(S^2\ot \sigma)(\Delta(h)(1\ot h))=S^2(h)\ot \sigma(h)$$
and if we apply $\varepsilon$ on the second factor in this equality and use that $S^2(h)=h$, we find that $\sigma(h)\varepsilon(\sigma(h))=h\varepsilon(\sigma(h))$. So, if $\varepsilon(\sigma(h))\neq 0$ we must have that $\sigma(h)=h$. In that case we of course must have $\varepsilon(\sigma(h))=1$ as $\varepsilon(h)=1$. This proves the first statement.
\snl
To prove the second statement, we use that $(h\ot 1)\Delta(h)=h\ot h$. We get
$$(S^2\ot \sigma)(h\ot h)=(S^2\ot \sigma)((h\ot 1)\Delta(h))=(S^2(h)\ot 1)\Delta(\sigma(h))$$
and because $S^2(h)=h$, we find by applying $\varepsilon$ on the second factor, that $h\sigma(h)=0$ because we assume that 
$\varepsilon(\sigma(h))=0$. Finally $\varphi(hah)=\varphi(ah\sigma(h))=0$.
\hfill$\square$
\einspr

Because $\sigma'(a)=\delta\sigma(a)\delta^{-1}$, we have $\varepsilon(\sigma'(h))=\varepsilon(\sigma(h))$ as $\varepsilon(\delta)=1$. If this  number is non-zero, we will not only have $\sigma(h)=h$ but also $\sigma'(h)=h$. This is shown in a similar way. But then, it also follows that $h$ and $\delta$ will commute. 
\snl
In the $^*$-algebra case and when the integrals are positive, we can only have that $\varphi(hah)=0$ for all $a$ when $h=0$. So, also here, in this case we have $\sigma(h)=h$. 
\snl
The case where $\varepsilon(\sigma(h))=0$ seems to be quite exceptional. We have no examples where this is true. And in this paper, we in fact need that $\varepsilon(\sigma(h))\neq 0$. 
\nl
Therefore, most of the time, {\it we will assume this extra condition} so that $\sigma(h)=h$. We will call $h$ a {\it regular} group-like idempotent. In the other case, we call it {\it exceptional}.
\nl
In the $^*$-algebra case, we know from [L-VD] that two algebraic quantum hypergroups can be constructed from a group-like projection $h$ in $A$. One is obtained by taking $hAh$ for the underlying algebra with coproduct $\Delta_0(a)=(h\ot h)\Delta(a)(h\ot h)$ for $a\in hAh$. The other case is dual to this. We use that the Fourier transform $k$ in $\widehat A$, defined by $k=\varphi(\,\cdot\,h)$, is again a group-like projection in the dual $\widehat A$ (provided $\varphi$ is normalized so that $\varphi(h)=1$).
\nl
The {\it first construction} is still possible when $h$ is a group-like projection in the multiplier algebra $M(A)$. We will formulate this result as a theorem.
\snl
Before we do this, we want to make a remark. In general, when $B$ is a subalgebra of an algebra $A$ with a non-degenerate product, it can happen that the product in $B$ is degenerate. When the product in $B$ is also non-degenerate, we can consider the multiplier algebra $M(B)$ and look for a possible relation with the original multiplier algebra $M(A)$. In general, not much can be said. Certainly, if $m\in M(A)$ and if it satisfies $mB\subseteq B$ and $Bm\subseteq B$, then $m$ defines a multiplier $n$ of $B$. The map $m\mapsto n$ is an algebra map from a subalgebra of $M(A)$ to $M(B)$. There is in general no reason for this map to be either injective or surjective.
\snl
For the subalgebra $A_0$ we consider in the next theorem, it can be shown that this map is surjective but it will not be injective. For the subalgebra $A_1$ we will consider later, this map will be both injective and surjective (see a remark following Proposition 1.12 below).

\inspr{1.3} Theorem \rm Let $(A,\Delta)$ be an algebraic quantum group and $h$ a regular group-like idempotent in $M(A)$. Let $A_0=hAh$ and define $\Delta_0$ on $A_0$ by 
$$\Delta_0(a)=(h\ot h)\Delta(a)(h\ot h)$$ 
when $a\in A_0$. Then $(A_0,\Delta_0)$ is an algebraic quantum hypergroup. The counit, the antipode and the integrals are obtained by simply taking the restrictions of these objects on $A$ to the subalgebra $A_0$.

\snl\bf Proof\rm:
First, we need to argue that the product in $A_0$ is still non-degenerate. Assume e.g.\ that we have an element $a\in A_0$ and that $ab=0$ for all $b\in A_0$. Because $A$ has local units, there exists an element $e\in A$ so that $ae=a$. Because $a\in A_0$, we have $ah=a$ and therefore we get also $aheh=a$. Now, $heh\in A_0$ and then we have $aheh=0$ by assumption. So $a=0$. Similarly on the other side. Hence the product in $A_0$ is non-degenerate and we can consider the multiplier algebras $M(A_0)$ and $M(A_0\ot A_0)$.
\snl
We clearly can define $\Delta_0(a)$ in $M(A\ot A)$ for any $a\in A$  by the formula in the formulation of the theorem. Now, assume that $a\in A$ and $b\in A_0$. Then we have
$$\Delta_0(a)(1\ot b)=(h\ot h)\Delta(a)(h\ot h)(1\ot b)=(h\ot h)\Delta(a)(1\ot b)(h\ot h)$$
because $hb=b=bh$. It follows that $\Delta_0(a)(1\ot b)\in A_0\ot A_0$. Similarly for the three other cases. Among other things, this implies that $\Delta_0(a)$ is defined as a multiplier in $M(A_0\ot A_0)$. And because there is no difficulty with coassociativity, we see that $\Delta_0$ defines a regular coproduct on $A_0$.
\snl
Because $\varepsilon(h)=1$, we see that the restriction of $\varepsilon$ to $A_0$ will be a counit for $\Delta_0$ on $A_0$.
\snl
We claim that the restriction of the left integral $\varphi$ of $A$ to the subalgebra $A_0$ is still a left integral. Indeed, because we have a group-like idempotent, we have 
$$\Delta_0(a)=(h\ot 1)\Delta(hah)(h\ot 1)=(h\ot 1)\Delta(a)(h\ot 1)$$
when $a\in A_0$ and when we apply $\varphi$ on the second leg of this equation, we get from the left invariance of $\varphi$ on $A$ that 
$$(\iota\ot\varphi)\Delta_0(a)=\varphi(a)h$$
for all $a\in A$. Because $h$ is the identity in $M(A_0)$, this gives left invariance of $\varphi$ also on $(A_0,\Delta_0)$.
\snl
The restriction of $\varphi$ to $A_0$ is non-trivial because we assume that the group-like idempotent is regular. Indeed, if $\varphi(hah)=0$ for all $a$, then $\varphi(ah)=0$ for all $h$ because $h\sigma(h)=h^2=h$. Because $\varphi$ is faithful on $A$, this would imply $h=0$. 
\snl
By a similar argument, we find that the restriction of the right integral $\psi$ will give the right integral on $(A_0,\Delta_0)$.
\snl
Finally, we need to show that there is an antipode relative to $\varphi$ on $A_0$. Because it is assumed that $S(h)=h$, we will also have that $S$ leaves $A_0$ globally invariant. And because
$$(\iota\ot \varphi)(\Delta_0(a)(1\ot b))=h((\iota\ot\varphi)(\Delta(a)(1\ot b))h$$
for $a,b\in A_0$  and similarly when the factor $b$ is on the other side, we get the necessary equation 
$$S((\iota\ot\varphi)(\Delta_0(a)(1\ot b)))=(\iota\ot\varphi)((1\ot a)\Delta_0(b))$$
whenever $a,b\in A_0$, simply because this formula is true when $\Delta_0$ is replaced by $\Delta$.
\snl
All this proves that the pair $(A_0,\Delta_0)$ is an algebraic quantum hypergroup.
\hfill$\square$
\einspr 

In the case of a $^*$-algebraic quantum group with a  regular group-like idempotent $h$ satisfying $h^*=h$, it is clear that the subalgebra $A_0$ is $^*$-invariant and that $\Delta_0$ is a $^*$-map so that we get a $^*$-algebraic quantum hypergroup. If moreover the integrals on $A$ are positive, then the same will be true for the integrals on $A_0$. 
\snl
It is not so difficult to obtain the other data. The modular automorphisms are also the restrictions of the modular automorphisms of the original algebraic quantum group (as $\sigma(h)=h$ and $\sigma'(h)=h$ so that they will leave the subalgebra globally invariant). And because the modular element $\delta$ of $A$  commutes with $h$, it follows that the modular element of the algebraic quantum hypergroup $(A_0,\Delta_0)$ is just the product $h\delta$.
\snl
Remark that in [L-VD], we found $h\delta=h$ but this is not necessarily true here. If e.g.\ $h=1$ and $\delta\neq 1$, this will not hold. In the $^*$-case, with positive integrals, we necessarily have that the scaling constant $\nu$ has to be $1$. If $h\in A$ and if $\varphi(h)\neq 0$, then it will also follow that $\nu=1$ as $S^2(h)=h$. But in general, this need not be true if $h$ is exceptional. 
\nl
Let us now look at the {\it second construction}.
\snl
In general, the Fourier transform $k$ of a group-like idempotent in $M(A)$, as defined above, will be an element in the dual space $A'$ (of all linear functionals on $A$), but it need not be an element even of $M(\widehat A)$. Indeed, take the trivial case where $h=1$. This is certainly a group-like projection, the Fourier transform will be $\varphi$ itself and this only belongs to $M(\widehat A)$ when $A$ is unital. Therefore, we can not proceed as in the case of a group-like projection in $A$ as in [L-VD] to construct the dual case.
\snl
What we will do in the remainder of this section is to start with a group-like projection $k$ in $\widehat A$ instead and use this and the ideas from [L-VD] to construct an associated algebraic quantum hypergroup. If $k$ turns out to be the Fourier transform of a group-like projection $h\in A$ as above, then this algebraic quantum hypergroup will be precisely the dual one mentioned before. See Example 2.3 Section 2.
\nl
\it Conditional expectations associated to a group-like idempotent in $M(\widehat A)$ \rm
\nl
In what follows, we will use the left and right $\widehat{A}$-bimodule structure on $A$. We know that these module structures on $A$ are unital in the sense that $\widehat{A} \blacktriangleright A = A$ and $A \blacktriangleleft \widehat{A} = A$.  As mentioned in the introduction already, these actions extend to $M(\widehat{A})$ in a natural way.
\snl
In the next propositions we prove that the existence of a group-like idempotent in $M(\widehat{A})$ is equivalent to the existence of specific linear maps from $A$ to itself.

\inspr{1.4} Proposition \rm 
Let $k$ be a group-like idempotent in $M(\widehat{A})$. We can define linear maps $E$, $E'$ from $A$ to itself by
$$E(a) = k \blacktriangleright a \qquad\qquad E'(a) = a \blacktriangleleft k.$$
We have that $E(E(a)) = E(a)$ for all $a\in A$ and that $E(E(a)a') = E(a) E(a')$ and also $ E(a E(a'))=E(a)E(a')$ for all $a,a'\in A$. So, $E$ is a conditional expectation on $A$. The same is true for $E'$.

\bf\snl Proof\rm:
Because of the remark preceding the formulation of the proposition, and because we assume that $k\in M(\widehat A)$, we have that $E$ and $E'$ are well-defined linear maps from $A$ to itself.
\snl
Because $k$ is assumed to be an idempotent, it follows immediately that $E(E(a))=E(a)$ and $E'(E'(a))=E'(a)$ for all $a\in A$.  
\snl
Now let $a,a'\in A$. Then we have 
$$E(E(a)a') = k \btr ((k \btr a) a') = \sum_{(k)}(k_{(1)}k \btr a)(k_{(2)}\btr a') = 
(k\btr a)(k\btr a').$$
and so also $E(E(a)a')=E(a) E(a')$. In the same way, we can show that
 $E(aE(a')) = E(a) E(a')$. So $E$ is a conditional expectation on $A$.  
\snl
The proof for the map $E'$ is similar. 
\hfill$\square$
\einspr

We now obtain some other properties of these two conditional expectations $E$ and $E'$ associated to the group-like idempotent $k\in M(\widehat A)$. 

\inspr{1.5} Proposition \rm 
With the assumptions of the previous proposition we have
\snl
i) $EE' = E'E$, \newline
ii) $E \circ S = S \circ E'$ and $E' \circ S = S \circ E$, \newline
iii) $(E\ot \iota)\Delta=(\iota\ot E')\Delta.$

\bf\snl Proof\rm:
i) Follows from the fact that $A$ is a $\widehat{A}$-bimodule so that left and right actions commute. To prove the first statement in ii), take $a\in A$. We have 
$$E(S(a)) = k \btr S(a) = S(a \btl k) = S(E'(a))$$
(where we have used that $S(k)=k$). Similarly for the second statement.
\snl
To prove the last equality we use the Sweedler notation. By definition we have 
$$E(a)=\sum_{(a)}a_{(1)}k(a_{(2)}) \qquad\qquad\text{and}\qquad\qquad E'(a)=\sum_{(a)}k(a_{(1)})a_{(2)}$$
for all $a$ and so
$$(E\ot \iota)\Delta(a)=(\iota\ot E')\Delta(a)=\sum_{(a)}a_{(1)}k(a_{(2)})a_{(3)}.$$
\vskip -0.7 cm\hfill$\square$
\einspr

The equation iii) can be given a meaning in different ways as can be seen from the proof. We just need to cover $a_{(1)}$ and $a_{(3)}$. This can be done e.g.\ by applying elements $f,g\in \widehat A$ on the first and the last factor respectively.
The equation is then read as $(fk)g=f(kg)$ in $\widehat A$. 
\snl
This form of the equation will be used in the proof of the following proposition where we obtain a converse of the above results.

\inspr{1.6} Proposition \rm
Suppose that $E$ and $E'$ are conditional expectations on $A$ satisfying the following properties:
\snl
i) $E \circ S = S \circ E'$ and $E' \circ S = S \circ E$, \newline
ii) $(E\ot \iota)\Delta=(\iota\ot E')\Delta.$
\snl
Then there is a group-like idempotent $k$ in $M(\widehat{A})$ such that $E$  and $E'$ are given by the left, respectively the right action of $k$ on $A$.

\snl\bf Proof\rm: 
Let $E$ and $E'$ be conditional expectations on $A$ as in the formulation of the proposition. It follows from the second condition that we can define  a multiplier $k$ of $\widehat A$ by the formulas
$$fk = f \circ E \qquad\qquad\text{and}\qquad\qquad kf = f\circ E'$$
where $f\in \widehat A$. Indeed, as we have seen in the remark following the previous result, we can rewrite this condition as $(fk)g=f(kg)$ in $\widehat A$ an this is precisely what we need in order to be able to define the multiplier $k$ of $\widehat A$.
\snl
From the very definition of $k$, we see that $E(a) = k \btr a$ and $E'(a) = a \btl k$ whenever $a\in A$. 
\snl  
To show that $k$ is an idempotent, take  $f \in \widehat{A}$. Then  we have $fk^2 = f \circ E^2 = f \circ E = fk$ and therefore  $k^2 = k$.
\snl
Next we show that $S(k)=k$. For all $a \in A$ and $f\in \widehat A$, we have, using the first condition, that
$$\align
(f(\widehat{S}(k)))(a) &=  (k \widehat{S}^{-1}(f))(S(a))= (\widehat{S}^{-1}(f))( E' (S(a)))\\
&= \widehat{S}^{-1}(f)(S(E(a))) = f(E(a))=  (fk)(a)
\endalign$$
and because this is true for all $a\in A$ and all $f\in \widehat A$, we see that $\widehat{S}(k) = k$.
\snl
Finally we prove $\widehat{\Delta}(k) (1 \otimes k) = k \otimes k$. The proof of the other formula is similar.
For all $f \in \widehat{A}$ and $a, b \in A$, we have
$$\align
	\langle (\widehat{\Delta}(fk)(1 \otimes k)), a \otimes b\rangle &=  \langle \widehat{\Delta}(fk), a \otimes E(b)\rangle\\
	&= \langle  fk,  aE(b)\rangle = \langle f,E(aE(b))\rangle.
\endalign$$
Because $E$ is a conditional expectation, we have $E(aE(b))=E(a)E(b)$ and so we get 
$$f(E(aE(b))=f(E(a)E(b))=\langle \widehat \Delta(f),E(a)\ot E(b)\rangle= \langle\widehat \Delta(f)(k\ot k),a\ot b\rangle.$$
As this is true for all $a,b\in A$ we get
$$\widehat{\Delta}(fk)(1 \otimes k)=\widehat \Delta(f)(k\ot k)$$
and because this is true for all $f\in \widehat{A}$, we find $\widehat{\Delta}(k)(1 \otimes k)=k\ot k.$
This completes the proof.
\hfill$\square$
\einspr

Remark that we automatically have that $E$ and $E'$ commute. In fact, it follows immediately from condition ii) in the previous proposition that $E$ commutes with the right action and that $E'$ commutes with the left action of any element in $\widehat A$.
\nl
Now recall the following definition (see e.g.\ Definition A.13 in [L-VD]). A subalgebra $B$ of $A$ is called {\it left invariant} if $\Delta(B)(A\ot 1)$ and $(A\ot 1)\Delta(B)$ are subsets of $A\ot B$. Similarly a right invariant subalgebra is defined. 

\inspr{1.7} Proposition \rm
The space $E(A)$ is a left-invariant subalgebra of $A$ and the space $E'(A)$ is a right invariant subalgebra of $A$.

\snl \bf Proof\rm:
First remark that $E(A)$ and $E'(A)$ are subalgebras because they are the range of a conditional expectation.
\snl
For $a, a' \in A$, we have 
$$\Delta(E(a)) (a' \otimes 1) = \Delta(k \btr a)(a' \otimes 1) = \sum a_{(1)} a' \otimes (k \btr a_{(2)}).$$ 
Therefore we have $\Delta(E(A))(A \otimes 1) \subseteq A \otimes E(A)$. The inclusion $(A \otimes 1)\Delta(E(A)) \subseteq A \otimes E(A)$ is proven in a similar way.
\snl
The right-invariance of the subalgebra $E'(A)$ is proven similarly.
\hfill$\square$
\einspr
It follows automatically that also $\Delta(E(A))(1 \otimes E(A)) \subseteq A \otimes E(A)$ and similarly for $E'(A)$. Indeed, we have $\Delta(E(A))(A \otimes E(A)) \subseteq A \otimes E(A)$ because $E(A)$ is an algebra and as $\Delta(E(A))(1 \otimes E(A)) \subseteq A \otimes A$, this can only happen when already $\Delta(E(A))(1 \otimes E(A)) \subseteq A \otimes E(A)$. We will need this form of invariance later.
\nl
We will now investigate the relation between the integrals and these conditional expectations on $A$. As before, $\varphi$ denotes a left integral on $A$ and $\psi$ a right integral.

\inspr{1.8} Proposition \rm 
Let $k$ be a group-like idempotent in $M(\widehat A)$ and let $E,E'$ be the associated conditional expectations on $A$. 
We always have $\varphi \circ E' = \varphi$ and $\psi \circ E=\psi$. If moreover the group-like idempotent is regular, we also have $\varphi \circ E = \varphi$ and $\psi\circ E'=\psi$. If on the other hand it is exceptional, we get $\varphi \circ E =0$ and $\psi\circ E'=0$.

\bf\snl Proof\rm:
For all $a \in A$, we have
$$\align \varphi(E(a))&= \varphi(k \btr a) = \varphi(a)\langle k,\delta \rangle\\
        \varphi(E'(a))&= \varphi(a \btl k) = \varphi(a)\langle k, 1 \rangle
\endalign$$
and because $\langle k,\delta \rangle=\widehat\varepsilon(\widehat\sigma(k))$ (see e.g.\ Proposition 4.1 in [De-VD]) and $\langle k, 1 \rangle=\widehat\varepsilon(k)$, all the results follow for the left integral. The proof for the right integral is completely similar (or follows if we use that the antipode converts $E$ to $E'$ and $\varphi$ to $\psi$).
\hfill$\square$
\einspr

Conversely, we see from the above that regularity of $k$ follows if $\varphi$ and $\psi$ are invariant under both $E$ and $E'$. This is another reason why it seems that regularity is very natural and that the other case is indeed exceptional.
\nl
Next, we look for the relations of $E$ and $E'$ with the modular element and the modular automorphisms. The results follow easily from the invariance of the integrals under these conditional expectations.  
\snl
We first make the following observation. We assume that $k$ is regular. For all $a,b\in A$, we get
$$\align \varphi(E(a)b)&=\varphi(E(E(a)b))=\varphi(E(a)E(b))\\
			\varphi(aE(b))&=\varphi(E(aE(b)))=\varphi(E(a)E(b))
\endalign$$
so that $\varphi(E(a)b)=\varphi(aE(b))$. Similar formulas hold for $E$ replaced by $E'$ and for $\varphi$ by $\psi$ as the two integrals are invariant with respect to the two conditional expectations. 

We now show that $E$ and $E'$ commute with $\sigma$ and $\sigma'$.

\inspr{1.9} Proposition \rm
Let $k$  be a regular group-like idempotent. Then
$$\align E(\sigma(a))&=\sigma(E(a)) \,\qquad \qquad \text{and} \qquad \qquad E(\sigma'(a))=\sigma'(E(a)) \\
         E'(\sigma(a))&=\sigma(E'(a)) \qquad \qquad \text{and} \qquad \qquad E'(\sigma'(a))=\sigma'(E'(a)) 
\endalign$$
for all $a\in A$.

\snl\bf Proof\rm:
Take $a,b\in A$. From the above observation we find
$$\varphi(aE(\sigma(b)))=\varphi(E(a)\sigma(b))=\varphi(bE(a))=\varphi(E(b)a)=\varphi(a\sigma(E(b)))$$
and by the faitfulness of $\varphi$ we find $E(\sigma(b))=\sigma(E(b))$. 
\snl
The three other formulas follow in the same way.
\hfill$\square$
\einspr

Next, we show that $E$ and $E'$ commute with multiplication with $\delta$.

\iinspr{1.10} Proposition \rm 
Let $k$ be a regular group-like idempotent. Then
$$\align
&E(a) \delta = E(a\delta) \qquad \qquad \text{and} \qquad \qquad E'(a)\delta = E'(a\delta) \\
&\delta E(a) = E(\delta a) \qquad \qquad \text{and} \qquad \qquad \delta E'(a) = E'(\delta a)
\endalign$$
for all $a \in A$.
\snl\bf Proof\rm:
We assume that $\psi$ and $\varphi$ are chosen so that $\psi=\varphi(\,\cdot\,\delta)$. Then for all $a,b\in A$, again using the observation above, we get
$$\varphi(aE(b\delta))=\varphi(E(a)b\delta)=\psi(E(a)b)=\psi(aE(b))=\varphi(aE(b)\delta)$$
and again from the faithfulness of $\varphi$ we find $E(b\delta)=E(b)\delta$.
\snl
The three other formulas are proven in the same way.
\hfill$\square$
\einspr

There are of course several other possible ways to obtain these relations, e.g.\ by using the definitions of $E$ and $E'$ as in Proposition 1.4. and it is instructive to verify this. We leave it as an exercise for the reader. It is also not so clear what is still possible if the group-like idempotent $k$ is exceptional.
\nl
Now, we will proceed as in Section 3 of [L-VD] and construct an algebraic quantum hypergroup from a pair of conditional expectations satisfying the previous properties. 
\nl
\it The associated algebraic quantum hypergroup and its dual \rm
\nl
So, fix an algebraic quantum group $(A,\Delta)$ and a  group-like idempotent $k$ in $M(\widehat A)$. We assume moreover that $k$ is regular so that $\widehat\sigma(k)=k$. Let $E, E'$ be the associated conditional expectations as in Proposition 1.4.
\snl
In Proposition 1.7 we have seen that the images $E(A)$ and $E'(A)$ are subalgebras of $A$. Then this is also true for the intersection of these two subspaces. This intersection is precisely the image of the composition $EE'$ because $E$ and $E'$ commute (cf.\ Proposition 1.5). 

\iinspr{1.11} Notation \rm
We denote the subalgebra $E(A)\cap E'(A)$ by $A_1$. So $A_1=EE'(A)=E'E(A)$.
\einspr

Remark that this subalgebra is invariant under the antipode because the antipode converts $E$ to $E'$ and vice-versa (see again Proposition 1.5). 
\snl
The algebra $A_1$ will be the underlying algebra of the algebraic quantum hypergroup we will construct. First we show that the algebra $A_1$ is still, in some sense, large enough (and therefore can not be completely trivial).

\iinspr{1.12} Proposition \rm
With $A$ and $A_1$ as above, we have
$$AA_1 = A_1A = A.$$

\snl\bf Proof\rm: 
Let $a_1, a_2, ... a_n$ be a finite number of elements in  $A$. Choose $b \in E(A)$ such that $\psi(b) = 1$.  This is possible because $\psi = \psi \circ E$. Consider the expressions $(1 \otimes a_i) \Delta(b)$ in $A \otimes A$.  By the existence of local units in a regular multiplier Hopf algebra, there is an element $p \in A$ such that  $(1 \otimes a_i) \Delta(b) (p \otimes 1) = (1 \otimes a_i) \Delta(b)$ for all $i$.  By evaluating $\psi \otimes \iota$ on both sides of this equation, we obtain $a_ix=a_i$ where $x= (\psi \otimes \iota)(\Delta(b) (p \otimes 1))$.  As $b \in E(A)$ and $\Delta(E(A))(A \otimes 1) \subseteq A \otimes E(A)$, we obtain that $x\in E(A)$. So, we find (right) local units for $A$ in $E(A)$.
\snl
Now let $a\in A$ and take again $b\in E(A)$ so that $\psi(b)=1$. Similarly as above, but now using local units for $A$ in $E(A)$, we find that $ax=a$ where $x= (\psi \otimes \iota)(\Delta(b) (p \otimes 1))$ with now $p\in E(A)$. Then 
$x = S(y)$ where $y = (\psi \otimes \iota)((b \otimes 1) \Delta(p))$ and for the same reason as above, we will have $y\in E(A)$ because now $p\in E(A)$.  Because $S(E(A))=E'(A)$, it follows that also $x\in E'(A)$. Therefore we have that $x \in E(A) \cap E'(A) = (EE') (A)$. We have obtained that $AA_1 = A$.  
\snl
The proof for $A_1A = A$ is similar. 
\hfill $\square$
\einspr

We see from the proof that we get a slightly stronger result. There exist local units for $A$ in $A_1$.
\snl
It is also a consequence of this result that the multiplications in the algebras $E(A)$, $E'(A)$ and $A_1$ are still non-degenerate. Indeed, if e.g.\ $a\in A_1$ and $ab=0$ for all $b\in A_1$, then because $A=A_1A$, we also get $ab=0$ for all $b\in A$. Therefore $a=0$. Similarly on the other side and for the algebras $E(A)$ and $E'(A)$.
\snl
As these subalgebras have non-degenerate products, we can consider their multiplier algebras. And as a further consequence of these results, we can consider e.g.\ $M(A_1)$ as sitting in $M(A)$. Indeed, given $m\in M(A_1)$ we can put $m(ab)=(ma)b$ for $a\in A_1$ and $b\in A$ and similarly on the other side. In fact, $M(A_1)$ is characterized as those elements $m\in M(A)$ satisfying $ma\in A_1$ and $am\in A_1$ for all $a\in A_1$. Similarly for the other algebras $E(A)$ and $E'(A)$.
\snl
Observe the difference with the situation of the subalgebra $hAh$, discussed earlier in this section.
\nl
We now consider the {\it coproduct} $\Delta$ on this subalgebra $A_1$. 
\snl
In general, we do not have $\Delta(A_1) \subseteq M(A_1 \otimes A_1)$.  Consider e.g.\ the multiplier Hopf algebra $A = K(G)$ of complex valued functions with a finite support on a group $G$. The dual $\widehat{A}$ is the ordinary group algebra $\Bbb C[G]$.  For a finite subgroup $H$ in $G$, we have the group-like idempotent $k$ in $\Bbb C[G]$, given as $k = \frac{1}{n} \sum_{r \in H} r$ where $n$ is the number of elements in $H$. For all $f \in K(G)$ and $s \in G$, we have 
$$(k \btr f) (s) = \frac{1}{n} \sum_{r \in H} f(sr) \quad \text{and} \quad (f \btl k)(s) = \frac{1}{n} \sum_{r \in H} f(rs).$$  Clearly $A_1 = (EE')(A) = k \btr A \btl k$ is given by the subalgebra of all functions in $K(G)$ which are constant on the double cosets of $H$ in $G$. For $f,g \in A$ and $s$, $t$ in $G$, we have $\Delta(f)(1\otimes g)(s,t) = f(st)g(t)$.  However when $f$ and $g$ are taken in $A_1$, we do not have that $\Delta(f) (1\otimes g)$ sits in $A_1 \otimes A_1$ because we can not expect that $f(srt)=f(st)$ for all $r\in H$.
\snl
In order to get a coproduct on $A_1$, we have to cut it down by the conditional expectations in an appropriate way. The situation is the same as in Section 3 of [L-VD] and similar as in the case of the algebra $A_0$ in Theorem 1.3. We use iii) of Proposition 1.5 for the following intermediate step in our construction.

\iinspr{1.13} Proposition \rm We can define a linear map $\Delta_1:A\to M(E(A) \otimes E'(A))$ by
$$\Delta_1(a)=(E\ot\iota)\Delta(a)=(\iota\ot E')\Delta(a).$$

\snl\bf Proof\rm:
We have explained already that the expressions on the right hand side of this formula have a meaning. 
\snl
Take $a\in A$. For any $b\in E'(A)$, we will have $\Delta_1(a)(1\ot b)=(\iota\ot E')(\Delta(a)(1\ot b))$ because $E'$ is a conditional expectation. It follows that $\Delta_1(a)(1\ot b)\in E(A)\ot E'(A)$. Similarly $(1\ot b)\Delta_1(a)\in E(A)\ot E'(A)$. For any  $c\in E(A)$, using that $E$ is a conditional expectation, we get $\Delta_1(a)(c\ot 1)\in E(A)\ot E'(A)$, as well as $(c\ot 1)\Delta_1(a)\in E(A)\ot E'(A)$. All this implies that $\Delta_1(a)\in M(E(A) \otimes E'(A))$.
\hfill$\square$
\einspr

From this result, we get easily that $A_1$ with the coproduct $\Delta_1$ has the structure of an algebraic quantum hypergroup. Observe that for an algebraic quantum hypergroup, it is not assumed that the coproduct is an algebra map (see the introduction and [De-VD]).

\iinspr{1.14} Proposition \rm The map $\Delta_1$ defines a regular comultiplication on the subalgebra $A_1$. The restriction of the counit $\varepsilon$ of $A$ to $A_1$ is a counit for $(A_1,\Delta_1)$.

\snl\bf Proof\rm:
First, we prove that $\Delta_1(a) (1 \otimes b)$ sits in $A_1 \otimes A_1$ for all $a,b\in A_1$. 
\snl
We have seen in the previous proposition that $\Delta_1(a)(1\ot b)\in E(A)\ot E'(A)$ because $b\in E'(A)$. Because $a\in E(A)$ and $E(A)$ is left invariant (as in Proposition 1.7), and as also $b\in E(A)$, we have $\Delta(a)(1\ot b)\in A\ot E(A)$ (see also the remark following Proposition 1.7. This gives us that also $\Delta_1(a)(1\ot b)\in E(A)\ot E(A)$. Finally, because $a\in E'(A)$ and because $E'(A)$ is right invariant, we find $\Delta_1(a)(1\ot b)\in E'(A)\ot E'(A)$ and so $\Delta_1(a)(1\ot b)\in E'(A)\ot E'(A)$. Therefore we get $\Delta_1(a) (1 \otimes b)$ sits in $A_1 \otimes A_1$.  
\snl
In a completely similar way, we can show the three other inclusions needed to have a regular coproduct.
\snl
To prove coassociativity, take $a,b,c \in A_1$. Then we have
$$\align (c \otimes 1 \otimes 1) (\Delta_1 \otimes \iota)(\Delta_1 (a) (1 \otimes b))
 &= (E \otimes \iota \otimes E') ((c \otimes 1 \otimes 1)(\Delta \otimes \iota) (\Delta(a) (1 \otimes b))) \\
 &= (E \otimes \iota \otimes E')((\iota \otimes \Delta)((c \otimes 1) \Delta(a))(1\otimes 1 \otimes b)) \\
 &= (\iota \otimes \Delta_1)((c \otimes 1)\Delta_1 (a)) (1 \otimes 1 \otimes b).
\endalign$$
\snl
Finally, we have for all $a\in A_1$ that 
$$(\varepsilon\ot\iota)\Delta_1(a)=(\varepsilon\ot\iota)(\iota\ot E')\Delta(a)=E'((\varepsilon\ot\iota)\Delta(a))=E'(a)=a$$
and similarly for the other equation.
\hfill $\square$
\einspr

Now we are ready to prove the main result of this paper.

\iinspr{1.15} Theorem \rm 
Let $(A,\Delta)$ denote an algebraic quantum group.  Suppose that $E$ and $E'$ are conditional expectations associated to a regular group-like idempotent $k$ in the multiplier algebra $M(\widehat A)$ of the dual $\widehat A$. The pair $(A_1,\Delta_1)$ with $A_1=E(A)\cap E'(A)$ and $\Delta_1$ on $A_1$ as defined by $\Delta_1=(E\ot \iota)\Delta$, is an algebraic quantum hypergroup. The integrals and the antipode on $(A_1, \Delta_1)$ are given by the restrictions to $A_1$ of the integrals and the antipode on $A$.
\snl
If $(A, \Delta)$ is a $^\ast$-algebraic quantum group and $k$ is self-adjoint in $M(\widehat{A})$, then $(A_1, \Delta_1)$ is a $^\ast$-algebraic quantum hypergroup.  

\snl\bf Proof\rm: 
We have seen already that $A_1$ is an algebra with a non-degenerate product and that $\Delta_1$ is a regular coproduct on $A_1$. Following [De-VD, Definition 1.10], we still have to prove the existence of a faithful left integral on $A_1$ and an antipode relative to this integral.
\snl
We consider the restriction of the left integral $\varphi$ on $A$ to the subalgebra $A_1$.  We still have that $\varphi_{|{A_1}}$is faithful on $A_1$.  Indeed, if e.g.\ $a \in A_1$ and $b\in A$, we have $\varphi(aEE'(b))=\varphi(ab)$ and so if $\varphi(ac)=0$ for all $c\in A_1$, we can take $c=EE'(b)$ and it will follow from the faithfulness of $\varphi$ on $A$ that  $a=0$.
\snl
Further, it is easy to see that $\varphi$ is left-invariant relative to $\Delta_1$.
\snl
Next, we consider the restriction of the antipode $S$ to the algebra $A_1$.  Observe again that $S(A_1) = A_1$ because $E \circ S = S \circ E'$ and $E' \circ S = S \circ E$.  Then, for all $a, b \in A_1$, we have 
$$\align S(\iota \otimes \varphi)(\Delta_1(a)(1\otimes b)) &= S(\iota \otimes \varphi)(\Delta(a)(1 \otimes b))\\
      &= (\iota \otimes \varphi)((1 \otimes a) \Delta(b))\\
      &= (\iota \otimes \varphi) ((1 \otimes a) \Delta_1 (b)).
\endalign$$
\snl
Finally, assume that $(A, \Delta)$ is a $^\ast$-algebraic quantum group and $k$ is self-adjoint in $M(\widehat{A})$. The maps $E$ and $E'$ on $A$ are $^\ast$-maps and  therefore $A_1$ is a $^\ast$-subalgebra of $A$ and $\Delta_1$ is a $^\ast$-map.
\hfill $\square$
\einspr

Because $(A_1, \Delta_1)$ is now shown to be an algebraic quantum hypergroup, we can use the results of Section 3 in [De-VD] to get the dual algebraic quantum hypergroup $(\widehat{A_1}, \widehat{\Delta_1})$.
\snl
The elements in $\widehat A_1$ are linear functionals on $A_1$, of the form $\varphi(\,\cdot\,a)$ where $a\in A_1$ because the left integral on $A_1$ is the restriction of the left integral $\varphi$ on $A$. This means that elements in $\widehat A_1$ are also restrictions to $A_1$ of elements in $\widehat A$.
\snl   
In the following proposition we prove that the algebraic quantum hypergroup $(\widehat{A_1}, \widehat{\Delta_1})$ is isomorphic with the algebraic quantum hypergroup $((\widehat A)_0,(\widehat \Delta)_0)$ where $(\widehat A)_0=k\widehat A k$ and $(\widehat \Delta)_0=(k\ot k)\widehat \Delta(\,\cdot\,)(k\ot k)$ (as in Theorem 1.3).

\iinspr{1.16} Theorem \rm Let $k$ be a regular group-like idempotent in the multiplier algebra $M(\widehat A)$ of the dual $\widehat A$ of an algebraic quantum group $A$. Let $E$ and $E'$ be the conditional expectations associated to  $k$. Let $(\widehat A_1\widehat\Delta_1)$ be the dual of $(A_1,\Delta_1)$ in the sense of algebraic quantum hypergroups.  There is a natural isomorphism of $(\widehat A_1,\widehat\Delta_1)$ with $((\widehat A)_0,(\widehat \Delta)_0)$. 

\snl\bf Proof\rm: 
Define the canonical map $\gamma$ from $\widehat A_1$ to $\widehat A$ by $\gamma(\omega)=\varphi(\,\cdot\,a)$ on $A$ if $\omega=\varphi(\,\cdot\,a)$ with $a\in A_1$. Remark that the left integral on $A_1$ is the restriction of the left integral on $A$ and we are using the same symbol for the two integrals here.
\snl
We first show that $\gamma$ is a bijection from $\widehat A_1$ to the subalgebra $(\widehat A)_0$. Let $a\in A$ and denote $\varphi(\,\cdot\,a)$ in $\widehat A$ by $\hat a$. Then, for all $x\in A$ we have 
$$\langle k\hat a k,x\rangle=\langle \hat a,EE'(x)\rangle=\varphi(EE'(x)a)=\varphi(xEE'(a))$$
and we see that $k\hat a k$ is the Fourier transform of $EE'(a)$. Because $A_1$ is defined as $EE'(A)$, it follows that $\gamma$ is indeed a bijection from $\widehat A_1$ to the subspace $k\widehat A k$ of $\widehat A$.
\snl
Next we show that $\gamma$ is an isomorphism of algebras. Take $\omega,\omega'\in k\widehat A k$ and $x\in A$. Then
$$\align \langle \omega\omega',x\rangle &= \langle \omega\ot \omega',\Delta(x)\rangle \\
&=\langle \omega k\ot \omega',\Delta(x)\rangle\\
&=\langle \omega \ot \omega', (E\ot \iota)\Delta(x) \rangle \\
&=\langle \omega \ot \omega', \Delta_1(x) \rangle.
\endalign$$
This precisely means that $\gamma$ is an algebra map.
\snl
Finally, we show that $\gamma$ converts the coproduct $\widehat\Delta_1$ on $\widehat A_1$ to the coproduct $(\widehat \Delta)_0$ on $(\widehat A)_0$. Take $\omega,\omega'\in k\widehat A k$ and $x,y\in A$. Then we have
$$\align \langle \widehat\Delta(\omega)(1\ot \omega'), x\ot y \rangle 
		& =\langle \omega \ot \omega', (x\ot 1)\Delta(y) \rangle \\
		& =\langle \omega \ot k\omega', (x\ot 1)\Delta(y) \rangle \\
		& =\langle \omega \ot \omega', (\iota\ot E')((x\ot 1)\Delta(y)) \rangle \\
		& =\langle \omega \ot \omega', (x\ot 1)\Delta_1(y) \rangle. 
\endalign$$
It follows from this that $\widehat\Delta(\gamma(\omega))(1\ot \gamma(\omega'))=(\gamma\ot\gamma)((\widehat\Delta)_0(\omega)(1\ot \omega'))$
for all $\omega,\omega'\in \widehat A_1$.
\snl
This completes the proof.
\hfill $\square$
\einspr

Because the counit and the antipode are unique in an algebraic quantum hypergroup, we must have that $\gamma$ converts the counit and the antipode on $\widehat A_1$ to the counit and the antipode on $k\widehat Ak$. This is easily verified as in the two cases, they are the restrictions of the these maps to the subalgebra.
\snl
One also verifies that the isomorphism $\gamma$ is a $^*$-isomorphism if we are in the case of a $^*$-algebraic quantum group with a self-adjoint group-like idempotent $k$. The $^*$-structure on $\widehat A_1$ is defined by the same formula that is also used to define the $^*$-structure on $\widehat A$ and the $^*$-structure on $k\widehat Ak$ is the restriction of the $^*$-structure on $\widehat A$ itself. This easily implies that indeed $\gamma$  will be a $^*$-isomorphism.

\nl\nl


\bf 2. Special cases and examples of algebraic quantum hypergroups \rm
\nl
In this section, we consider {\it some special cases} and {\it examples}. We use the notations and assumptions of the previous section.
\snl
For some of the examples, we consider known cases and show how they fit into our theory. For some of them, the group-like idempotent is actually in the algebra itself. For others, we do have that the idempotent is not in the algebra but really in the multiplier algebra.
\snl
A few of our examples involve groups. Therefore, let us fix the following conventions and notations.

\inspr{2.1} Notation \rm
Assume that $G$ is a group. We use $K(G)$ for the algebra of complex functions with finite support on $G$ with its usual structure of multiplier Hopf $^*$-algebra. For any $p\in G$, we denote by $\delta_p$ the function in $K(G)$ that is $1$ on $p$ and $0$ everywhere else. We use $\Bbb C[G]$ for the group algebra of $G$ over $\Bbb C$ with its usual structure of Hopf $^*$-algebra. We use $p\mapsto \lambda_p$ for the imbedding of $G$ in $\Bbb C[G]$. Remark that $(K(G), \Bbb C[G])$ is a dual pair and that the pairing is given by the formula $\langle f,\lambda_s\rangle=f(s)$ for all $s\in G$.
\einspr
First we look at the {\it motivating example} as mentioned already in the abstract.

\inspr{2.2} Example \rm
Let $G$ be a group and $A=K(G)$. Let $H$ be a finite subgroup of $G$ and define $k = \frac{1}{n} \sum_{r \in H} \lambda_r$ where $n$ is the number of elements in $H$. Then $k$ is a group-like idempotent in $\widehat A$.  The conditional expectations $E$ and $E'$ on $A$ are given as follows.  For all $f \in A$ and $p\in G$, we have
$$
E(f) (p) = \frac{1}{n} \sum_{r \in H} f(pr)\qquad\qquad
E'(f)(p) = \frac{1}{n} \sum_{r \in H} f(rp).
$$
The algebra $A_1 = EE'(A)$ is given by the complex valued functions with finite support in $G$  that  are constant on the double cosets of $H$ in $G$. The comultiplication $\Delta_1$ on $A_1$ is given  by the formula
$$
\Delta_1(f)(p,q) = \frac{1}{n} \sum_{r \in H} f(prq).
$$
when $f \in A_1$ and $p, q \in G$. For the dual algebraic quantum hypergroup  $(\widehat{A_1}, \widehat{\Delta_1})$ we have $\widehat{A_1}= k (\Bbb C[G]) k$ and  $\widehat{\Delta}_1 (k \lambda_p k) = k \lambda_p k \otimes k \lambda_p k$ for all $p\in G$. This means that
$$\widehat{\Delta}_1 (k \lambda_p k)=\frac{1}{n^2}\sum_{r,s\in H} \lambda_{rps}\ot \frac{1}{n^2}\sum_{r,s\in H} \lambda_{rps}$$
for all $p\in G$.
\hfill$\square$
\einspr

Next, we take the case studied already in [L-VD]. 

\inspr{2.3} Example \rm
Consider any algebraic quantum group $(A,\Delta)$. Assume that $k$ is in fact in $\widehat A$ (and that it is regular). Then it will be the Fourier transform $\varphi(\,\cdot\,h)$ of a (regular) group-like idempotent $h\in A$ (assuming that $\varphi$ is normalized so that $\varphi(h)=1$). The conditional expectations $E$ and $E'$, as defined in Proposition 1.4, are given by
$$E(a)=(\iota\ot\varphi)(\Delta(a)(1\ot h))
		\qquad\quad\text{and}\qquad\quad
  E'(a)=(\varphi\ot\iota)(\Delta(a)(h\ot 1))$$
for $a\in A$. These are precisely the conditional expectations as obtained and studied already in [L-VD]. In that paper, the algebra $A_1$ is denoted by $C_1$ and also the restricted coproduct $\Delta_1$ has been considered there. It is given by
$$\Delta_1(a)=(\iota\ot\varphi\ot\iota)(\Delta^2(a)(1\ot h\ot 1))$$
for all $a\in A_1$. The resulting algebraic quantum hypergroup is of discrete type (see Theorem 3.17 in [L-VD]). 
The left cointegral is the group-like idempotent $h$.
\snl
Also the duality we obtain  in Theorem 1.16 is  already present in [L-VD] for this special case. The dual is $((\widehat A)_0,(\widehat\Delta)_0)$ where $(\widehat A)_0=k\widehat A k$ and $(\widehat\Delta)_0(a)=(k\ot k)\widehat\Delta(a)(k\ot k)$. Now we obtain an algebra with an identity, namely $k$, and the resulting dual algebraic quantum hypergroup is of compact type.
\snl
In fact, the above is only completely correct if $A$ is a $^*$-algebraic quantum group with positive integrals and if the idempotents are self-adjoint. However, it is easily seen that the $^*$-operation and the self-adjointness of $h$  and $k$ are not essential (as long as they are not exceptional - a case that can not occur in the theory of [L-VD]).
\hfill$\square$
\einspr

Of course, in Example 2.2, we have a special case of the situation in Example 2.3.
\snl
One may wonder if we can get examples in the group case (as in Example 2.2) where the idempotent is not in the algebra but in the multiplier algebra. Clearly, as in Example 2.2, the dual $\widehat A$ is an algebra with identity, we can not get this here. However, if we look at the dual, we can obtain such a case.

\inspr{2.4} Example \rm
Again take a group $G$ and now let $A$ be the group algebra $\Bbb C[G]$. The dual $\widehat A$ is $K(G)$. When $H$ is any subgroup of $G$, we let $k$ be the function on $G$ that is $1$ on $H$ and $0$ everywhere else. This will be a group-like idempotent in $M(\widehat A)$ and if $H$ is not finite, this will not belong to $\widehat A$. The corresponding conditional expectations $E$ and $E'$ coincide on $\Bbb C[G]$ and are given by
$$E(\sum_{p\in G} a(p) \lambda_p)=\sum_{r\in H} a(r) \lambda_r$$
when $a$ is any complex function on $G$ with finite support. Clearly in this case $A_1$ is nothing else but the group algebra of $H$ while the dual is the function algebra $K(H)$. We get a dual pair of algebraic quantum groups.
\hfill$\square$
\einspr

Of course, the example above is not very interesting. It is a special case of the following. Again, we consider the notations and assumptions as in Theorem 1.15 and Theorem 1.16 of the previous section.

\inspr{2.5} Proposition \rm 
If $k$ is assumed to be central in $M(\widehat A)$, then the algebraic quantum hypergroups $(A_1,\Delta_1)$ and $((\widehat A)_0,(\widehat\Delta)_0)$ are actually algebraic quantum groups. They are still dual to each other (now in the sense of algebraic quantum groups).

\snl \bf Proof\rm:
As $k$ is central we have $\omega k=k\omega$ for all $\omega\in \widehat A$. This means that $\omega(E(a))=\omega(E'(a))$ for all such $\omega$ and all $a\in A$. Then it follows that $E=E'$. Therefore $A_1=E(A)=E'(A)$. Now consider the restricted coproduct $\Delta_1$. It is given by the formula in Proposition 1.13. If  $a,b\in A_1$, we get 
$$\Delta_1(a)(1\ot b)=(E\ot\iota)(\Delta(a)(1\ot b))=(E'\ot\iota)(\Delta(a)(1\ot b))$$
and because $E'(A)$ is right invariant, we see that that $(E'\ot\iota)(\Delta(a)(1\ot b))=\Delta(a)(1\ot b)$. Therefore, for all $a,b\in A_1$, we have $\Delta_1(a)(1\ot b)=\Delta(a)(1\ot b)$. From this, it follows easily that $\Delta_1$ is still an algebra map. This proves that $(A_1,\Delta_1)$ is an algebraic quantum group.
\snl
The dual case is even simpler. Because $k$ is central in $M(\widehat A)$, we see that $(\widehat A)_0=\widehat A k$ and that $(\widehat\Delta)_0(b)=\widehat \Delta(b)(k\ot k)$ for all $b\in \widehat A$. Again we see that $(\widehat\Delta)_0$ is still an algebra map in this case. So, also $((\widehat A)_0,(\widehat\Delta)_0)$ is an algebraic quantum group.
\snl
That these two algebraic quantum groups are still each others dual in the sense of algebraic quantum groups follows trivially from the fact that the dual in the theory of algebraic quantum hypergroups is constructed in the same way as for algebraic quantum groups and of course that the integrals on the restrictions above are simply the restrictions of the integrals.
\snl
In fact, one of the arguments above would be sufficient to prove the result by the last statement about duality.	
\hfill$\square$
\einspr  

Next, we consider {\it the example with double cosets} of a pair of finite-dimensional Hopf algebras as given by Vainerman in [V].

\inspr{2.6} Example \rm 
Let $A$ and $B$ be finite-dimensional Hopf algebras and let $\pi : A \rightarrow B$ be an epimorphism of Hopf algebras.  We suppose that $A$ and $B$ are co-Frobenius (i.e.\ that they have integrals).  Recall that the antipode of a co-Frobenius Hopf algebra is always bijective.  Further we suppose that $\varphi_B (1) = 1$ where $\varphi_B$ is the left integral on $B$. This implies that $B$ is unimodular in the sense that the left integral and the right integral coincide. 
\snl
We take $k=\varphi_B \circ \pi$ in $A'$.  Notice that $\widehat{A} = A'$ because $A$ is finite-dimensional.
\snl
It is easily seen that $k$ is a group-like idempotent in $A'$ and $\langle k, \delta_A\rangle = \varphi_B (\pi(\delta_A))$.  Because $\pi$ is a surjective homomorphism of Hopf algebras, we must have $\pi(\delta_A)=\delta_B$ and because $B$ is unimodular, we have $\delta_B=1$. Therefore $\langle k, \delta_A\rangle  = 1$ so that $k$ is a regular group-like idempotent.
\snl
By the use of Proposition 1.4, we have the following conditional expectations on $A$:
$$\align
E(a) &= k \btr a = (\iota \otimes \varphi_B)((\iota \otimes \pi) \Delta(a))\\
E'(a) &= a \btl k = (\varphi_B \otimes \iota)((\pi \otimes \iota) \Delta(a))
\endalign$$
for all $a \in A$.
\snl
The algebra $E(A)$ is given by the invariants of the right coaction $\gamma_r : A \rightarrow A \otimes B$ where $\gamma_r (a) = (\iota \otimes \pi) \Delta(a)$.  Similarly, $E'(A)$ is given by the invariants of the left coaction $\gamma_\ell : A \rightarrow B \otimes A$ where $\gamma_\ell(a) = (\pi \otimes \iota)\Delta(a)$.  By the use of Theorem 1.15, we have the algebraic quantum hypergroup $(A_1, \Delta_1)$ where $A_1 = (EE')(A)$ and $\Delta_1(a) = (E\otimes \iota) \Delta(a)=(\iota \otimes E')\Delta(a) = \sum \varphi_B (\pi(a_{(2)})) a_{(1)} \otimes a_{(3)}$. 
\snl
By Theorem 1.16, we have that $(\widehat{A_1},\widehat{\Delta_1})$ is given as $k \widehat{A}k$ and $\widehat{\Delta}_1 (\omega) = (k \otimes k) \widehat{\Delta}(\omega)(k \otimes k)$ for all $\omega \in k \widehat{A}k$.  It is easily seen that $k \widehat{A}k$ is the set of all linear functionals $f$ on $A$ such that $f = f \circ EE'$.
\snl
Both algebraic quantum hypergroups are finite-dimensional. They are of compact type and of discrete type. The algebras are unital and there are cointegrals. The unit in $A_1$ is simply the unit of $A$ as $EE'(1)=1$ and the unit in $k \widehat{A}k$ is of course $k$ itself. As $k \widehat{A} k = \widehat{A_1} = \{\varphi(\cdot\ (EE')(a))\mid a \in A\}$, we know that $k = \varphi_B \circ \pi = \varphi_A(\cdot\ h)$ for a unique element $h \in (EE')(A)$. Because $k$ is the unit in $\widehat{A_1}$, we have that $h$ is a left cointegral for $(A_1, \Delta_1)$ such that $\varphi_A(h) = 1$. The left cointegral in $k \widehat{A}k$ is $\varphi_A$.
\hfill$\square$
\einspr

Next, we consider an example studied in [DC]. We again use it to illustrate the procedure to construct a dual pair of algebraic quantum hypergroups as in Section 1. It is a case where the underlying algebraic quantum group is actually a $^*$-algebraic quantum group $(A,\Delta)$ with positive integrals and where in general, the group-like idempotent $k$ sits in $M(\widehat A)$ and not in $\widehat A$ itself. So, it does not fit into the framework of [L-VD].

\inspr{2.7} Example \rm
The starting point is a $^*$-algebraic quantum group $(C,\Delta_C)$ with a positive left integral $\varphi_C$. 
\snl
It is shown in [DC-VD] that $C$ is spanned by the eigenvectors of $S^2$  and that the eigenvalues are strictly positive real numbers. Moreover, for each eigenvalue $t$ of $S^2$, there is a map $P_t$ from $C$ to $C$, projecting onto the space of eigenvectors with eigenvalue $t$. It is uniquely determined by the requirement that it is a self-adjoint projection in the sense that $\varphi_C(P_t(a)^*b)=\varphi_C(a^*P_t(b))$ for all $a,b\in C$. The same result is true for the dual $\widehat C$ and the set of eigenvalues of $S^2$ is the same for $C$ and its dual $\widehat C$.
\snl
The dual $\widehat C$ will be denoted by $D$ in what follows and the data associated will be indexed with $C$ and $D$ respectively. However, we will use $P_t$ for the projection on the space of eigenvectors with eigenvalue $t$ for both $C$ and $D$.
\snl
Denote by $G$ the (multiplicative) subgroup of $\Bbb R^+\setminus \{0\}$ generated by the eigenvalues of $S^2$. We consider the dual pair $(K(G),\Bbb C[G])$ as in Notation 2.1. 
\snl
A new $^*$-algebraic quantum group $(A,\Delta_A)$ is constructed as follows. For $A$, we take the $^*$-algebra $A$, defined as $\Bbb C[G]\ot C$ with the usual tensor product $^*$-algebra structure. The coproduct $\Delta_A$ however is a twisted coproduct. It is given by the formula
$$\Delta_A(\lambda_p\ot c)=\sum_{t,(c)}\lambda_p\ot P_t(c_{(1)})\ot \lambda_{t^{-1}p} \ot c_{(2)}$$
where $p\in G$,  where $c\in C$ and the Sweedler notation is used for $\Delta_C$ and finally, where the sum over $t$ runs over all eigenvalues of $S_C^2$. The counit $\varepsilon_A$ and the antipode $S_A$ on $(A,\Delta)$ are given by
$$\align \varepsilon_A(\lambda_p\ot c)	&=\delta(p,1)\varepsilon_C(c)\\
	 S_A(\lambda_p\ot c)&=\lambda_{p^{-1}t}\ot S_C(c)
\endalign$$
for all $p\in G$  and $c\in C$ satisfying $S_C^2(c)=tc$. We used the Kronecker delta. A left integral $\varphi_A$ on $A$ is given by 
$$\varphi_A(\lambda_p\ot c)=\varphi_C(c)$$
for all $p\in G$ and $c\in C$.
\snl
We can also describe the dual $(\widehat A,\widehat \Delta)$ of $(A,\Delta)$. We will use $(B,\Delta_B)$ to denote this dual. The algebra $B$ is the space $K(G)\ot D$ with a twisted product. It is given by 
$$(\delta_p \otimes d)(\delta_q \otimes d') = \delta_p \delta_{t q} \otimes dd'$$
when $p, q \in G$ and when $d\in  D$ satisfes $S_D^2(d) = t d$ for some non-zero positive real number $t$. The involution is given by $(\delta_p\otimes d)^\ast = \delta_{t^{-1}p} \otimes d^\ast$ when again $S_D^2(d)=t d$. The coproduct is the usual tensor coproduct. The counit and the antipode on $(B,\Delta_B)$ are given by 
$$\align
\varepsilon_B(\delta_p \otimes d)& = \delta(p,1)\varepsilon_D(d)\\
S_B(\delta_p \otimes d) &= \delta_{p^{-1}t} \otimes S_D(d)\\
\endalign$$
when $p \in G$ and  $S_D^2(d) = t d$. A left integral $\varphi_B$ on $B$ is given by 
$$\varphi_B(\delta_p\ot d)=\varphi_D(d)$$
for all $p\in G$ and $d\in D$.
\snl
Now we consider $k$, defined as $\delta_1\ot 1$ in $M(B)$ (where we use $1$ for the number $1$ in the first place and $1$ for the identity in $M(D)$ in the second place). One easily verifies that it is indeed a multiplier and that
$$\align k(\delta_p\ot d)&=\delta(p,1)(\delta_1\ot d) \\
         (\delta_p\ot d)k&=\delta(p,t)(\delta_t\ot d) \qquad\text{if} \qquad S_D^2(d)=t d
\endalign$$
where $p\in G$ and $d\in D$. Moreover, it is easy to check that it is a group-like idempotent as in Definition 1.1.
\snl
It follows  from the above formulas that the algebra $kBk$ is
$$\{\delta_1\ot d \mid d\in D \ \text{with}\ S_D^2(d)=d \}$$
and the reduced coproduct on $kB k$ is given by the formula
$$\Delta_0(\delta_1\ot d)=\sum_{(d)}\delta_1 \ot P_1(d_{(1)}) \ot \delta_1 \ot d_{(2)}$$
for $d\in D$ satisfying $S_D^2(d)=d$. We see that the resulting $^*$-algebraic quantum hypergroup is isomorphic with $(D_0,\Delta_0)$ where $D_0=\{d\in D \mid S_D^2(d)=d \}$ and where $\Delta_0$ is defined on $D_0$ by
$$\Delta_0(d)=(P_1\ot\iota)\Delta_D(d).$$
\snl
Next, we look at the dual. The conditional expectations $E$ and $E'$ associated with $k$ on $A$ are given by the formulas
$$E(\lambda_p\ot c)= \lambda_p\ot P_p(c)
		\qquad\quad\text{and}\qquad\quad
			E'(\lambda_p\ot c)=\delta(p,1)(\lambda_1\ot c)$$
where $p\in G$ and $c\in C$. One verifies that 
$$EE'(\lambda_p\ot c)=\delta(p,1)(\lambda_1\ot P_1(c))$$
for all $p\in G$ and $c\in C$. The reduced coproduct is given by
$$\Delta_1(\lambda_1\ot c) =\sum_{(c)}\lambda_1 \ot P_1(c_{(1)}) \ot \lambda_1 \ot c_{(2)}$$
when $S^2_C(c)=c$. We see that now we get a $^*$-algebraic quantum hypergroup, isomorphic with $(C_0,\Delta_0)$  where $C_0=\{c\in C \mid S_C^2(c)=c \}$ and where $\Delta_0$ is defined on $C_0$ by
$$\Delta_0(c)=(P_1\ot\iota)\Delta_C(c).$$
\hfill$\square$
\einspr

It is interesting to remark that, although the two ways of constructing an algebraic quantum hypergroup in the example above are essentially different, they do give essentially the same type of algebraic quantum hypergroup.
\snl
In any case, the example illustrates very well our theory.
\nl
Finally, we come to our last example, constructed from a matched pair $(G_1,G_2)$ of groups. We refer to [M1, M2] for the original theory of bicrossproducts and to [De-VD-W] for the extension of these results to multiplier Hopf algebras. We refer to the notations in 2.1.

\inspr{2.8} Example \rm 
Let $(G_1, G_2)$ be a matched pair of (arbitrary) groups in the sense of [M1, Example 3.14]. See also Definition 6.2.10 in [M2].
\snl
 Elements in $G_1$ will be denoted by letters $r,s, \dots$ and elements of $G_2$ by letters $u,v,\dots$. The identity in each of these groups will be denoted by $e$. The left action of $G_1$ on $G_2$ is denoted as $r\tr u$ for $r\in G_1$ and $u\in G_2$ while the right action of $G_2$ on $G_1$ is denoted as $r\tl u$ when $r\in G_1$ and $u\in G_2$. By assumption, these actions satisfy the following relations.
For all  $s, t \in G_1$ and $u, v \in G_2$ we have
$$\align 
s \triangleright (uv) &= (s\triangleright u)((s  \triangleleft u)\triangleright v) \\ (st) \triangleleft u &= (s \triangleleft (t \triangleright u))(t \triangleleft u).
\endalign$$
We also have 
$s \triangleright e = e$ and $e \triangleleft u = e$ for every $s\in G_1$ and $u\in G_2$.
\nl
One can consider {\it two bicrossproducts}.  
\nl
i) On the one hand, there is the  {\it right-left bicrossproduct} $(A, \Delta_A)$. Here $A$ is the smash product $\Bbb C[G_1]\# K(G_2)$, constructed with the left action of $G_1$ on $G_2$. The product satisfies 
$$(\lambda_r \# \delta_u) (\lambda_s \# \delta_v) = \delta(u,s \triangleright v) (\lambda_{rs} \# \delta_v)$$
for $r, s \in G_1$ and $u, v \in G_2$ (where we also use $\delta$ for the Kronecker delta function). The  $^\ast$-operation is given by $(\lambda_r \# \delta_u)^\ast = \lambda_{r^{-1}} \# \delta_{r \triangleright u}$. The coproduct $\Delta_A$ is given by
$$\Delta_A (\lambda_r \# \delta_u) = \sum\limits_{w \in G_2} (\lambda_r \# \delta_w) \otimes (\lambda_{r \triangleleft w} \# \delta_{w^{-1}u})$$
whenever $r\in G_1$ and $u\in G_2$. The counit and the antipode on $(A,\Delta_A)$ are given by 
$$\varepsilon_A (\lambda_r\# \delta_u)=  \delta(u,e) 
	\qquad\quad \text{and} \qquad\quad 
		S_A(\lambda_r \# \delta_u) = \lambda_{(r \triangleleft u)^{-1}} \# \delta_{(r\triangleright u)^{-1}}.
$$
\snl
A left integral $\varphi_A$ is given by $\varphi_A (\lambda_r \# \delta_u) = \delta(r,e)$ whenever $r\in G_1$ and $u\in G_2$. It is also right invariant. 
\nl
ii) 
The dual $(\widehat{A},\widehat{\Delta_A})$, in the sense of algebraic quantum groups, of the algebraic quantum group $(A,\Delta)$ as described in item i) above, is given as the {\it left-right bicrossproduct}. We denote it by $(B,\Delta_B)$ in what follows.
\snl
The algebra $B$ is the smash product $K(G_1) \# \Bbb C[G_2]$, now constructed with the right action of $G_2$ on $G_1$. For $r, s \in G_1$ and $u,v \in G_2$, we have
$$(\delta_r \# \lambda_u) (\delta_s \# \lambda_v) = \delta(r \triangleleft u, s) (\delta_r \# \lambda_{u v})$$
for $r,s\in G_1$ and $u,v\in G_2$. The involution is given by $(\delta_r \# \lambda_u)^\ast  = \delta_{r \triangleleft u} \# \lambda_{u^{-1}}$. The coproduct is given by
$$
\Delta_B (\delta_r \# \lambda_u) = \sum\limits_{t \in G_1} (\delta_{rt^{-1}} \# \lambda_{t \triangleright u}) \otimes (\delta_t \# \lambda_u)$$
when $r \in G_1$ and $u \in G_2$. The counit and the antipode satisfy 
$$\varepsilon_B(\delta_r \# \lambda_u)= \delta(e,r) 
	\qquad\quad\text{and} \qquad\quad
		S_B(\delta_r\# \lambda_u) = \delta_{(r\triangleleft u)^{-1}} \# \lambda_{(r \triangleright u)^{-1}}.
$$
A left integral $\varphi_B$  is given by the formula $\varphi_B(\delta_r \# \lambda_u) = \delta(u,e)$. Again here, the left integral is also right invariant. 
\nl
For details about these constructions and formulas, we refer to [De-VD-W, Example 1.8, Example 2.11 and Example 3.17]. See also [M2, Example 6.2.11 and Example 6.2.12].
\nl
iii) Further, we now {\it assume that $G_2$ is a finite group}. Then $(A, \Delta_A)$ is an ordinary Hopf algebra. However, the dual $(\widehat A,\widehat{\Delta_A})$, denoted as $(B,\Delta_B)$, is now a genuine multiplier Hopf algebra if the other group $G_1$ is infinite. Only when the two groups are finite, we get two Hopf algebras that are then necessarily finite-dimensional.
\snl
We consider the element $k$ in $M(B)$ defined as $k = 1 \# k_0$ where $k_0$ is the element in $\Bbb C[G_2]$ defined as $k_0= \frac{1}{n} \sum_{u \in G_2} \lambda_u$  where $n$ is the number of elements in $G_2$. The element $k_0$ is actually a left (and a right) cointegral, normalized so that $\varepsilon_B(k_0)=1$. 
\snl
We have 
$$\align k(\delta_r \# \lambda_u) &= \frac{1}{n} \sum_{v \in G_2} \delta_{r \triangleleft v} \# \lambda_{v^{-1} u}\\
	(\delta_r \# \lambda_u)k &= \frac{1}{n} \sum_{v \in G_2}\delta_r \#  \lambda_v = \delta_r \# k_0
\endalign$$ 
where $r\in G_1$ and $u\in G_2$. We easily see that $k^2 = k=k^*$ and $S_{B}(k) = k$. Moreover, a straightforward calculation shows that $\Delta_{B} (k) (1 \otimes k) = k \otimes k$ as well as $\Delta_{B} (k) (k\otimes 1)=k\ot k$ (where now $1$ is the identity in $M(B)$). Therefore, $k$ is a (self-adjoint) group-like idempotent in $M(B)$, in the sense of Definition 1.1. In the general case, that is when $G_1$ is infinite, it does not belong to $B$ itself and so the case is not covered in [L-VD]. Also observe that $k$ is central if and only if the right action of $G_2 $ on $G_1$ is trivial. As the modular element $\delta_A$ in $(A,\Delta_A)$ is given by the identity, we have $\langle k,\delta_A\rangle = \varepsilon_B (k) = 1$ so that this idempotent is a regular one.  
\nl
iv) We can now {\it apply our theory}. We first look at the algebra $B_0$, defined as $kBk$ and the coproduct $\Delta_0$ on $B_0$, defined as in Theorem 1.3.
\snl
A straightforward calculation gives that for any $r\in G_1$ and $u\in G_2$ we have
$$k(\delta_r \# \lambda_u)k=\frac{1}{n} \sum_{v \in G_2} (\delta_{r \triangleleft v} \# \lambda_{v^{-1} u})k=P(\delta_r)\# k_0$$
where we define 
$$P(\delta_r)=
\frac{1}{n}\sum_{v \in G_2}\delta_{r \triangleleft v}$$
for every $r\in G_1$.
We see that $B_0$, defined as $kBk$, consists of elements $f\# k_0$ where $f$ is a function in $K(G_1)$, invariant under the action of $G_2$. It is in fact isomorphic with the subalgebra $K_{\text{inv}}(G_1)$ of $K(G_1)$ of such functions. We will denote it by $D$.
\snl
For the coproduct $\Delta_0$, as defined in Theorem 1.3, we find for all $r\in G_1$ and $u\in G_2$,
$$\align (k\ot k)\Delta_B(\delta_r\#\lambda_u)(k\ot k)
	&=\sum\limits_{t \in G_1} k(\delta_{rt^{-1}} \# \lambda_{t \triangleright u})k\otimes 	
		k(\delta_t \# \lambda_u)k\\
	&=\sum\limits_{t \in G_1} (P(\delta_{rt^{-1}}) \# k_0) \otimes 	
		(P(\delta_t) \# k_0)
\endalign$$
and we see that 
$$\align  
\Delta_0(P(\delta_r)\# k_0)
&=\Delta_0(k(\delta_r\#\lambda_e)k)\\
&=\Delta_0(\delta_r\#\lambda_e)\
=\sum\limits_{t \in G_1} (P(\delta_{rt^{-1}}) \# k_0) \otimes (P(\delta_t) \# k_0)
\endalign$$
for all $r\in G_1$. For the counit $\varepsilon_0$ and the antipode $S_0$, obtained by restricting the count $\varepsilon_B$ and the antipode $S_B$, we find
$$\varepsilon_0(P(\delta_r)\# k_0)=\delta(r,e)
\qquad\quad\text{and}\qquad\quad
S_0(P(\delta_r)\# k_0)=P(\delta_{r^{-1}})\# k_0$$
for all $r\in G_1$. To obtain the last formula, we use
$$S_0(P(\delta_r)\# k_0)=S_B(k(\delta_r\# \lambda_e)k)=kS_B(\delta_r\#\lambda_e)k=k(\delta_{r^{-1}}\# \lambda_e)k=P(\delta_{r^{-1}})\# k_0.$$
\snl
We see that the resulting quantum hypergroup is isomorphic with the algebra $D$, defined as the space $K_{\text{inv}}(G_1)$ of functions in $K(G_1)$, invariant under the action of the group $G_2$ and with the coproduct $\Delta_D$ given by
$$\Delta_D(P(\delta_r))=\sum\limits_{t \in G_1} P(\delta_{rt^{-1}})\otimes P(\delta_t)
=(P\ot P)\Delta(\delta_r)\tag"(2.1)"$$
for all $r\in G_1$. Here the counit is given by $\varepsilon(P(\delta_r))=\delta(r,e)$, the antipode by $S(P(\delta_r))=P(\delta_{r^{-1}})$ and the left integral $\varphi$ is given by $\varphi(P(\delta_r))=1$ for all $r$.
\nl
v) We now consider the dual case. The conditional expectations $E$ and $E'$ on $(A, \Delta_A)$, defined in Proposition 1.4, are given by
$$\align
E(\lambda_r \# \delta_u) &= k\blacktriangleright (\lambda_r \# \delta_u) = \frac{1}{n} (\lambda_r \# 1)\\
E' (\lambda_r \# \delta_u) &= (\lambda_r \# \delta_u) \blacktriangleleft k =  \frac{1}{n} \sum_{w \in G_2} (\lambda_{r \triangleleft w} \# \delta_{w^{-1}u})
\endalign$$
for all $r\in G_1$ and $u\in G_2$. We see that 
$$(EE') (\lambda_r \# \delta_u) = (E'E) (\lambda_r \# \delta_u) =\frac{1}{n} Q(\lambda_r) \# 1$$
for $r\in G_1$ and $u\in G_2$, where now $Q(\lambda_r)=\frac{1}{n} \sum_{w \in G_2} \lambda_{r \triangleleft w}$. The map $Q$ is a projection map on the group algebra $\Bbb C[G_1]$ and its range consists of elements invariant under the right action of $G_2$. So, the algebra $A_1$, defined as $(EE')(A)$, consists of elements $a_0 \# 1$ where $a_0$ is in the group algebra $\Bbb C[G_1]$ such that $a_0\tl u=a_0$ for all $u\in G_2$.
\snl
For $\Delta_1$, as defined in Proposition 1.14, we find the following. For any $r\in G_1$ and $u\in G_2$ we have
$$\align
	\Delta_1( \lambda_r \# \delta_u)
         &= (E\ot \iota)\Delta_A(\lambda_r \# \delta_u) \\
			&= \frac{1}{n} \sum_{w\in G_2} (\lambda_r \# 1) \ot (\lambda_{r\tl w}\# \delta_{w^{-1}u}).
\endalign$$
Therefore 
$$\align \Delta_1( \lambda_r \# 1)
        			&= \frac{1}{n} \sum_{w\in G_2} (\lambda_r \# 1) \ot (\lambda_{r\tl w}\# 1)\\
					&=  (\lambda_r \# 1) \ot (Q(\lambda_r) \# 1)
\endalign$$
for all $r\in G_1$. It follows from this that
$$\Delta_1(Q(\lambda_r)\# 1)=(Q(\lambda_r)\# 1)\ot (Q(\lambda_r)\# 1)$$
for all $r\in G_1$. Notice that we get a coabelian algebraic quantum hypergroup whose underlying algebra has an identity (because already $A$ was an ordinary Hopf algebra).
\snl
The different data for this quantum hypergroup are, as in all the other cases, again obtained by restricting the original data to the subalgebra. We have for the counit $\varepsilon_1$ that $\varepsilon_1(Q(\lambda_r)\# 1)=1$ and for the antipode $S_1$ that
$$S_1(Q(\lambda_r)\# 1) = Q(\lambda_{r^{-1}}) \# 1$$
for all $r\in G_1$.
\snl
We see that in this case, the algebra $A_1$ is isomorphic with the subalgebra $C$ of $\Bbb C[G_1]$ of elements invariant under the action of $G_2$ and the coproduct $\Delta_C$ on this isomorphic subalgebra is given by
$$\Delta_C(Q( \lambda_r))
        			=  Q(\lambda_r) \ot Q(\lambda_r) \tag"(2.2)"$$
for all $r\in G_1$.
\hfill$\square$\einspr

Notice how the duality between $B_0$ and $A_1$, as we have it from our general result in Theorem 1.16, translates very nicely to the isomorphic algebras $D$ and $C$. In the one case, we have the subalgebra of $K(G_1)$ of functions invariant under the action of $G_2$ while in the other case, we have the subalgebra of $\Bbb C[G_1]$, again with functions invariant under this action. The pairing between these algebras $C$ and $D$ is obtained by restricting the original pairing between $K(G_1)$ and $\Bbb C[G_1]$. The coproducts $\Delta_C$ and $\Delta_D$ induced on these subalgebras in this restricted pairing are in the two cases obtained by the natural projection maps on these subalgebras, coming from the restriction of the pairings. This yields precisely the formulas (2.1) and (2.2) above.
\snl
One should compare the results in this example with those in Example 2.7. In the two cases, we get isomorphic quantum hypergroups that can be defined directly, but then they will not fit in our general construction method.
\snl
The above example will also work if we consider a general bicrossproduct $A\# B$ where $A$ is any algebraic quantum group and $B$ a finite-dimensional Hopf algebra. And surely more interesting examples can be obtained using the ideas above.
\nl
\nl


\bf References \rm
\nl
{\bf [DC]} K.\ De Commer : {\it Algebraic quantum hypergroups imbedded in algebraic quantum groups}. Preprint University Tor Vergata (Rome) (2009).

{\bf [DC-VD]} K.\ De Commer \& A.\ Van Daele: {\it Multiplier Hopf algebras imbedded in locally compact quantum groups}. Preprint K.U.\ Leuven (2006), Arxiv math.OA/0611872v2, to appear in the Rocky Mountain Journal of Mathematics.

{\bf [De-VD]} L.\ Delvaux \& A.\ Van Daele: {\it Algebraic quantum hypergroups}. Preprint University of Hasselt and K.U.\ Leuven (2006), ArXiv math.RA/0606466.

{\bf [De-VD-W]} L.\ Delvaux, A.\ Van Daele and S.H.\ Wang, {\it Bicrossproducts of multiplier Hopf algebras}. Preprint University of Hasselt, K.U.\ Leuven and Nanjing University (2009), Arxiv math.RA/09032974. 

{\bf [Dr-VD]} B.\ Drabant \& A.\ Van Daele: {\it Pairing and Quantum double of multiplier Hopf algebras}. Algebras and Representation Theory {\bf 4} (2001), 109-132.

{\bf [Dr-VD-Z]} B.\ Drabant, A.\ Van Daele \& Y.\ Zhang: {\it Actions of multiplier Hopf algebras}. Comm.\ Algebra {\bf 27} (1999), 4117-4172. 

{\bf [K]} A.\ A.\ Kalyuzhnyi: {\it Conditional expectations on quantum groups and new examples of quantum hypergroups}. Methods of Funct.\ Anal.\ Topol.\ {\bf 7} (2001), 49-68.

{\bf [L-VD]} M.B.\ Landstad \& A.\ Van Daele): {\it Compact and discrete subgroups of algebraic quantum groups}. Preprint University of Trondheim and University of Leuven (2007), Arxiv math.OA/0702458v1.

{\bf [M1]} S.\ Majid: {\it Physics for algebraists: non-commutative and non-cocommutative Hopf algebras by a bicrossproduct construction}. Journal of Algebra {\bf 130} (1990), 17-64.

{\bf [M2]} S.\ Majid: {\it Foundations of quantum group theory}. Cambridge University Press (1995).

{\bf [V]} L.I.\ Vainerman: {\it Gel'fand pair associated with the quantum groups of motions of the plane and $q$-Bessel functions}. Reports on Mathematical Physics {\bf 35} (1995), 303-326.

{\bf [VD1]} A.\ van Daele: {\it Multiplier Hopf algebras}, Trans.\ Am.\ Math.\ Soc. {\bf 342}(2)(1994), 917-932.

{\bf [VD2]} A.\ Van Daele: {\it An Algebraic framework for group duality}, Adv.\ in Math.\ {\bf 140} (1998), 323-366.


\end